\documentclass[a4paper,10pt,reqno]{article}
\usepackage{amsfonts}
\usepackage{amsmath, amsfonts, amssymb, amsthm}
\numberwithin{equation}{section}

\begin{document}

\title{Energy conservation of weak solutions and Onsager$'$s conjecture for the two-component Camassa-Holm system}

\author{Jingjing Liu \footnote{e-mail:
jingjing830306@163.com } and Zhaoyu Sun \footnote{e-mail: sunzhaoyu0623@163.com} \\
School of Mathematics and Information Science, \\
Zhengzhou University of Light Industry,\\
Zhengzhou 450002, Henan, People$'$s Republic of China}
\date{}
\maketitle

\begin{abstract}
In this paper, we investigate energy conservation for weak solutions to the two-component Camassa-Holm system. We derive a local energy balance equation in the sense of distributions with an explicit defect term. As an application, we determine that the Onsager exponent is 1, and as a consequence, we also obtain the uniqueness of weak solutions.\\

\end{abstract}

\section{Introduction}
In this paper, we consider the classical model arising in shallow water wave theory, namely the periodic two-component Camassa-Holm (2CH) system :
\begin{equation}
 \ \ \ \ \ \ \ \  \ \ \ \ \left\{\begin{array}{ll}
m_{t}+2u_xm+um_{x}+\sigma\rho\rho_{x}=0,&t > 0,\,x\in \mathbb{R},\\
 \rho_{t}+(\rho u)_x=0, &t > 0,\,x\in \mathbb{R},\\
u(t,x+1)=u(t,x), & t \geq 0, x\in \mathbb{R},\\
\rho(t,x+1)=\rho(t,x), & t\geq 0, x\in \mathbb{R},\end{array}\right. \\
\end{equation}
where $m=u-u_{xx}$, and $\sigma=\pm1.$  This system was derived by Constantin and Ivanov in the context of shallow water
theory\cite{C-I}, where $u(t,x)$ denotes the horizontal fluid velocity and $\rho(t,x)$ represents the horizontal deviation of the
surface from equilibrium. For physical validity, they impose that $u(t,x)\rightarrow 0$ and $\rho(t,x)\rightarrow 1$
as $|x|\rightarrow \infty$ at any instant $t$. In this article, since we are dealing with a periodic situation, there is no need to
perform a transformation like $\overline{\rho}=\rho-1.$ Furthermore, according to \cite{C-I}, when $\sigma=-1$ , it indicates that
the gravitational acceleration is upward, while when $\sigma=1$, it indicates that the gravitational acceleration is downward.
We let $\sigma=1$ in this paper.

For $\sigma=-1$, system (1.1) was derived earlier within mathematical theory as a generalization of the Camassa-Holm equation
 \cite{R.Camassa,Constantin 0,Constantin 03,EY,GY,YYG} (the
Camassa-Holm equation can be obtained via the obvious reduction $\rho=0$ in (1.1)). This extension stems from the deformation theory
of semisimple bihamiltonian structures of hydrodynamic type. Liu and Zhang \cite{L-Z} founded its basic framework and generalized
the Camassa-Holm hierarchy to two components by deforming compatible Poisson brackets. On this basis, Chen, Liu and Zhang \cite{C-L-Z}
 formalized the 2CH system, established its reciprocal transformation to the first negative flow of the AKNS hierarchy, and obtained
 peakon and multi-kink solutions. Meanwhile, Falqui\cite{F} independently investigated the 2CH system via Lie algebra Hamiltonian
 structures, yielding its Lax representations, N-peakon solutions consistent with the classic Camassa-Holm equation, and non-standard
 algebro-geometric solutions.

After system (1.1) was derived, it received widespread attention and yielded a wealth of results, see \cite{C-I,E-L-Y,G,GY1,GY2,GL1,GL2,LZ}.
These works established local well-posedness, derived precise
blow-up scenario, proved that the existence of
strong solutions which blow up in finite time, as well as the global existence of strong solutions \cite{C-I,E-L-Y,G,GY1,GL1,GL2}. In \cite{GY2}
and \cite{LZ}, the authors investigated the existence and uniqueness of weak solutions for system (1.1), respectively. In the study of classical
solutions, the conservation of energy of the solutions, that is
\begin{equation*}
E(t)=\int_{\mathbb{S}}(u^{2}(t,x)+u_{x}^{2}(t,x)+\rho^{2}(t,x))dx=\int_{\mathbb{S}}(u_{0}^{2}(x)+u_{0x}^{2}(x)+\rho_{0}^{2}(x))dx,
\end{equation*}
here $\mathbb{S}=\mathbb{R}/\mathbb{Z}\simeq [0,1)$ stands for the unit circle, $u_{0}(x)=u(0,x),$ $\rho_{0}(x)=\rho(0,x)$ are
initial value, plays a very important role. Then, after the blow up of the classical solutions, does the energy of the weak
solutions studied still remain conserved? This is a question worth considering. Generally speaking, there are two situations.
One is that energy dissipates, that is, $ E(t) < E(T)$ for any $t>T.$ The other situation is that the energy can still remain
conserved. Particularly, the author obtained uniqueness within the category of conserved weak solutions in \cite{LZ}.

Onsager$'$s conjecture states that for weak solutions of the 3D incompressible Euler equations, kinetic energy is conserved
if the spatial H$\ddot{o}$lder regularity exponent is greater than $\frac{1}{3}$, while anomalous energy dissipation can occur when
the exponent is less than $\frac{1}{3}$, see \cite{CWT,O}. There is a large body of work on the Onsager$'$s conjecture, including \cite{BGSTW,BT,CCFS,DE,DR,FW,I,MNY}.
Specifically, Boutros and Titi in \cite{BT} studied an analogue of Onsager$'$s conjecture for the 3D Euler-$\alpha$ equations (the inviscid Camassa-Holm equations).
They proved that the equation of local energy balance holds in the sense of distributions
with a precise defect term provided the weak solution $u\in L^{3}([0,T]; W^{1,3}(\mathbb{S}^{3})).$  Furthermore,
the energy is conserved as long as $u\in L^{3}([0,T]; B_{3,\infty}^{\beta}(\mathbb{S}^{3}))$ with $\beta>1$. This work represents the first time that an Onsager-type
analysis has been applied to the field of shallow water wave equations.

In \cite{WLL}, the authors studied the energy dissipation mechanism for weak solutions of various Camassa-Holm type equations
by establishing a local energy balance equation in the sense of distributions with a precise defect term. They obtained energy
conservation of the weak solutions under specific conditions and determined that their Onsager index is 1.  Inspired by this work, this paper
follows the method in \cite{WLL} to investigate energy conservation of weak solutions to the 2CH system (1.1). To the best of
our knowledge, this problem has not been addressed in the existing literature.

In contrast to single-component Camassa-Holm type equations, the two-component coupled 2CH system exhibits more intricate energy structures and dissipation behaviors. The main novelty of this paper lies in that we generalize the local energy balance method in the distributional sense to the periodic two-component shallow water wave model for the first time. We rigorously derive the explicit energy dissipation defect term for weak solutions of the 2CH system, establish sufficient conditions for energy conservation and the uniqueness of weak solutions, and determine its critical Onsager index. Moreover, the analytical framework constructed in this work can be extended to other multi-component integrable shallow water wave systems, which provides an effective tool for studying energy conservation laws of coupled nonlinear dispersive wave equations.

The outline of the remaining paper is as follows. In Section 2, we provide several key lemmas and outline the main conclusions of this work. Section 3 is devoted to the rigorous proof of Theorem 2.1, which details the construction of these terms for weak solutions to the periodic 2CH. Building on these constructions, Section 4 rigorously proves a corollary and a theorem from Section 2, establishing the regularity criteria for energy conservation and weak solution uniqueness.\\

\section{Preparations and main results}
\newtheorem{theorem2}{Theorem}[section]
\newtheorem{lemma2}{Lemma}[section]
\newtheorem{corollary2}{Corollary}[section]
\newtheorem {remark}{Remark}[section]
\newtheorem{definition2}{Definition}[section]
To facilitate the derivation of the energy conservation laws and the construction of the energy dissipation mechanisms, this section collects several preliminary analytical tools used throughout the paper.  Subsequently, we state the main results of this paper.

\begin{lemma2}
[see \cite{LF}]
Given a smooth cutoff function $\varphi(x)= C_0 e^{-\frac{1}{1-|x|^2}}$ for $|x| < 1$ and $\varphi(x)= 0$ for $|x| \geq 1 $, where the normalization constant $C_0 > 0$ satisfies $\int_{\mathbb{R}} \varphi(x) \, dx = 1$. For $\varepsilon > 0$, we define the rescaled mollifier $\varphi_\varepsilon(x) = \frac{1}{\varepsilon} \varphi\left(\frac{x}{\varepsilon}\right)$, we obtain $|\varphi_\varepsilon'(x)| = \frac{1}{\varepsilon^2} \left|\varphi'(\frac{x}{\varepsilon})\right|$. For any function $f\in L^1_{loc}(\Omega)$, its mollified version is defined as\\
$$
f^{\varepsilon}(x)=(f*\varphi_\varepsilon)(x)=\int_{\mathbb{R}}f(x-y)
\varphi_\varepsilon(y)dy,
$$
and $x\in \Omega_\varepsilon$,where
$$
\Omega_\varepsilon =\{x\in\Omega:dist(x,\partial \Omega)>\varepsilon\}.
$$
\end{lemma2}
\begin{lemma2}[see \cite{S}]
Suppose $1\leq p\leq \infty$ and $T \in ( {{L}^{p}( {\mathbb{R}}^{n}) ,{L}^{p}( {\mathbb{R}}^{n}) })$. Let $\widehat{u}$ be the multiplier
corresponding to $T$ and suppose that $\widehat{u}$ is continuous at each point of the lattice $\Lambda$. Set $\lambda (m) = \widehat{u}(m)$ for $m \in \Lambda$ , Then there exists a unique periodized operator $\widetilde{T}$ such
that $\widetilde{T} \in  \left( {{L}^{p}\left( {\mathbb{S}}^{n}\right) ,{L}^{p}\left( {\mathbb{S}}^{n}\right) }\right)$ and $\| \widetilde{T}\|  \leq  \| T\| .$\\
\end{lemma2}

Then we focus on system (1.1) with $\sigma=1$, that is
\begin{equation}
 \ \ \ \ \ \ \ \  \ \ \ \ \left\{\begin{array}{ll}
m_{t}+2u_xm+um_{x}+\rho\rho_{x}=0,&t > 0,x\in \mathbb{R},\\
 \rho_{t}+(\rho u)_x=0, &t > 0,x\in \mathbb{R},\\
u(t,x+1)=u(t,x), & t \geq 0,x\in \mathbb{R},\\
\rho(t,x+1)=\rho(t,x), & t\geq 0, x\in \mathbb{R},\\
u(0,x)=u_{0}(x), & x\in \mathbb{R},\\
\rho(0,x)=\rho_{0}(x), & x\in \mathbb{R}, \end{array}\right. \\
\end{equation}
and establish a local K\'{a}rm\'{a}n-Howarth-Monin type identity for weak solution to the system. Having established these preliminaries, we proceed to state the main theorems.

\begin{definition2}
Given $(u_{0},\rho_{0})\in H^{1}(\mathbb{S})\times L^{2}(\mathbb{S}).$ If $(u,\rho)\in L^{\infty}([0,T], H^{1}(\mathbb{S})\times L^{2}(\mathbb{S}))$ satisfies
system (2.1) in the sense of distribution, then $(u,\rho)$ is called a weak solution to system (2.1).
\end{definition2}

\begin{theorem2}
Assume $(u,\rho)$ is the weak solution corresponding to the initial data $(u_{0},\rho_{0})\in H^{1}(\mathbb{S})\times L^{2}(\mathbb{S})$ such that
$u\in L^\infty([0,T]; H^1(\mathbb{S}))$ and $u_x\in L^{2+}([0,T]; L^{2+}(\mathbb{S}))$, $\rho \in L^\infty([0,T]; L^2(\mathbb{S}))$ and $\rho_x \in L^{2+}([0,T]; L^{2+}(\mathbb{S}))$, then the following local K\'arm\'an--Howarth--Monin type equation holds in the sense of distributions
\begin{equation}
\begin{aligned}
&(u u^\varepsilon)_t + (u_x u_x^\varepsilon)_t + (\rho \rho^\varepsilon)_t
+ \frac{1}{2} (u^2 u^\varepsilon)_x
- [(u^2)^\varepsilon u]_x
+ [(u u_x) u_x^\varepsilon]_x\\
&+ [(\rho u) \rho^\varepsilon]_x+ (\Pi u^\varepsilon)_x + (\Pi^\varepsilon u)_x
+ \frac{1}{2} (u^3)^\varepsilon_x
+ \frac{1}{2} (u u_x^2)^\varepsilon_x\\
&- \frac{1}{2} [u (u_x^2)^\varepsilon]_x
+ \frac{1}{2} (u \rho^2)^\varepsilon_x
- \frac{1}{2} [u (\rho^2)^\varepsilon]_x \\
&= D^{2CH}_\varepsilon(u, \rho),
\end{aligned}
\end{equation}
with $\Pi=(1 - \partial_x^2)^{-1}\left(u^2 + \frac12 u_x^2 + \frac12 \rho^2\right)$,
where the K\'arm\'an-Howarth-Monin type relation and defined by \\
\begin{equation}
\begin{aligned}
D^{2CH}_\varepsilon(u, \rho)
&= -\frac{1}{2} \int_\mathbb{S} \varphi'_\varepsilon(l) [\delta u(l)]^3 dl
   -\frac{1}{2} \int_\mathbb{S} \varphi'_\varepsilon(l) \delta u(l) [\delta u_x(l)]^2 dl \\
&\quad -\frac{1}{2} \int_\mathbb{S} \varphi'_\varepsilon(l) \delta u(l) \left[\delta \rho(l)\right]^2 dl.
\end{aligned}
\end{equation}
Then
$${D}^{2CH}( {u,\rho}):=\lim\limits_{\varepsilon\rightarrow0}D_{\varepsilon}^{2CH} = {( {u}^{2}) }_{t} + (u^2_x)_t+ (\rho^2)_t + (u u^2_x)_x+{( {\rho}^{2}u) }_{x} + {2}{( {\Pi}{u}) }_{x}$$
in the sense of distributions.
\end{theorem2}

\begin{corollary2}
Assume that $u$, $\rho$ satisfy \\
\begin{equation}
          \int_{\mathbb{S}}|{\delta{u} (l)  }|^{3}{dx} \leq  {C}_{1}(t) |l| \sigma_1(l),
\end{equation}
\begin{equation}
      \int_{\mathbb{S}} |\delta{u}(l)| [{\delta{u}_{x}(l)}]^{2}{dx} \leq  {C}_{2}(t) | l| \sigma_2(l)
\end{equation}
and\\
\begin{equation}
\int_{\mathbb{S}}|{\delta{u}}( l)|[ { \delta\rho( l)   }]^{2}{dx} \leq  {C}_{3}( t) | l| \sigma_3( l),
\end{equation}
where $C_i(t)$ are integrable functions on $[0,T]$ and $\sigma_i(l)$ are bounded functions on some neighborhood of the origin for $i=1,2,3$. If in addition, $\sigma_i(l)\rightarrow 0 $ as $l\rightarrow 0$, then the energy $\|(u,\rho)\|_{{H_1(\mathbb{S})}\times{L^2(\mathbb{S})}}$ is conserved.
\end{corollary2}

\begin{theorem2}
 Suppose that $(u,\rho)$ is a weak solution of $(2.1)$.\ Then the energy is conserved as long as $(u_x,\rho)\in L^p([0,T];L^q(\mathbb{S})\times L^q(\mathbb{S}))$ with $\frac{1}{p}+\frac{2}{q}=1$ for $q\geq3$.
 Moreover, this weak solution is unique.
\end{theorem2}

\begin{remark}
\upshape
The regularity condition in Theorem 2.2 corresponds to an Onsager exponent of 1. This implies that the total energy of system (2.1) is conserved provided the weak solution possesses an Onsager index of exactly 1, thereby establishing an Onsager-type criterion tailored to the specific structure of (2.1).
\end{remark}

\section{Proof of Theorem 2.1 }
\newtheorem {remark3}{Remark}[section]
\newtheorem{theorem3}{Theorem}[section]
\newtheorem{lemma3}{Lemma}[section]
\newtheorem{definition3}{Definition}[section]
\newtheorem{claim3}{Claim}[section]
\newenvironment{proofof}[1]{%
  \par\noindent
  \textbf{Proof of  #1}%
  \quad
}{\hfill$\square$\par}

As a preliminary step toward the proof of the main result Theorem $2.1$, system $(2.1)$ is first reformulated in conservation law form.

Let $G(x) \mathrel{\text{:=}} \frac{\cosh ( {x-[x]- 1/2}) }{2\sinh \left( {1/2}\right) }$, $x\in \mathbb{R}$.\ Then ${(1-{\partial }_{x}^{2}) }^{-1}f = G*f$ for all $f\in {L}^{2}( \mathbb{S})$
and $G * m = u$.\ Here, we denote by $*$ the convolution.\ By a direct calculation, one can
rewrite system $(2.1)$ as follows:
\begin{equation}
 \left\{\begin{array}{ll}
{u}_{t} + u{u}_{x} + {\partial }_{x}G * ( {{u}^{2} + \frac{1}{2}{u}_{x}^{2} + \frac{1}{2}{\rho }^{2}})  = 0,&t > 0,\,x\in \mathbb{R},\\
 {\rho}_{t} + u{\rho}_{x} + {u}_{x}\rho = 0, &t > 0,\,x\in \mathbb{R},\\
u(t,x+1)=u(t,x), & t \geq 0, x\in \mathbb{R},\\
\rho(t,x+1)=\rho(t,x), & t\geq 0, x\in \mathbb{R},\\
u(0,x)=u_{0}(x), & x\in \mathbb{R},\\
\rho(0,x)=\rho_{0}(x), & x\in \mathbb{R}.\end{array}\right. \\
\end{equation}
By defining the non-local variable
$$G*(u^2 + \frac{1}{2} u_x^2 + \frac12 \rho^2) = \Pi = (1 - \partial_x^2)^{-1}(u^2 + \frac12 u_x^2 + \frac{1}{2} \rho^2),$$
$(3.1)_1$ and $(3.1)_2$  can be transformed into an equivalent formulation involving only first-order spatial dervatives of $\Pi$,
\begin{equation}
 \left\{\begin{array}{ll}
{u}_{t} + u{u}_{x} +\Pi_x = 0,&t > 0,\,x\in \mathbb{R},\\
 (1-\partial_{x}^2)\Pi = {{u}^{2} + \frac{1}{2}{u}_{x}^{2} + \frac{1}{2}{\rho}^{2}}, &t > 0,\,x\in \mathbb{R},\\
{\rho}_{t} + u{\rho}_{x} + {u}_{x}\rho  = 0, & t > 0, x\in \mathbb{R},\end{array}\right. \\
\end{equation}
By the second equation of $(3.2)$, we derive the elliptic relation,
\begin{equation}
\Pi_{xx}=\Pi-{u}^{2}-\frac{1}{2}{u}_{x}^{2}-\frac{1}{2}{\rho }^{2}.
\end{equation}
The equivalent system laying the foundation for deriving energy dissipation terms and establishing energy conservation results, next comes the proof of Theorem 2.1.
\vspace{1em}

\begin{proofof}{Theorem 2.1} In view of the regularity of the solutions, we need to perform the operations in the sense of distributions. For
$\psi\in C_{0}^{\infty}([0,T]\times\mathbb{S}),$ by $(3.2)_1$, we have
$$\langle u_{t}+uu_{x}+\Pi_{x}, \psi_{x}\rangle=0,$$
it follows that
$$\langle u_{tx}+u_{x}^{2}+{u} u_{xx}+\Pi_{xx}, \psi\rangle=0,$$
that is
$$u_{tx}+u_{x}^{2}+u u_{xx}+\Pi_{xx}=0$$ in the sense of distributions. The subsequent calculations are carried out in the sense of distributions.
For convenience, we omit the explicit writing of the test functions.

Substituting the expression for $\Pi_{xx}$ given in $(3.3)$ into this identity and combining it with $(3.2)_1$ and $(3.2)_3$, we obtain \\
\begin{equation}
\left\{\begin{array}{ll}
{u}_{t} + u{u}_{x} + {\Pi}_{x} = 0,\\
 {u}_{tx} + u{u}_{xx} + \frac{1}{2}{u}_{x}^{2} = {u}^{2} + \frac{1}{2}{\rho }^{2} - \Pi,\\
\rho_{t}+(\rho u)_{x}=0.\end{array}\right. \\
\end{equation}
Regularizing this system, we get
\begin{equation}
\left\{\begin{array}{ll}
{u}_{t}^{\varepsilon } + {( u{u}_{x}) }^{\varepsilon } + \Pi_{x}^{\varepsilon } = 0,\\
 {u}_{tx}^{\varepsilon} + {( u   {u}_{xx}) }^{\varepsilon} + \frac{1}{2}{( {u}_{x}^{2}) }^{\varepsilon} = {( {u}^{2}) }^{\varepsilon} + \frac{1}{2}{( {\rho }^{2}) }^{\varepsilon} - {\Pi}^{\varepsilon},\\
{\rho }_{t}^{\varepsilon} + {( \rho u) }_{x}^{\varepsilon} = 0.\end{array}\right. \\
\end{equation}
Multiplying $(3.4)_1$ and $(3.5)_1$ by $u^\varepsilon$ and $u$ respectively, we have
\[\left\{  \begin{array}{l} {u}_{t}   {u}^{\varepsilon} + {u}{u}_{x}   {u}^{\varepsilon} + {\Pi}_{x}  {u}^\varepsilon= 0,
\\  u   {u}_{t}^{\varepsilon} + {( u   {u}_{x}) }^{\varepsilon} u+ {\Pi}_{x}^{\varepsilon}u = 0 . \end{array}\right.\]
Multiplying $(3.4)_2$ and $(3.5)_2$ by $u_{x}^\varepsilon$ and $u_{x}$ respectively, we have
\[\left\{  \begin{array}{l} {u}_{xt}   {u}_{x}^{\varepsilon } + {u}{u}_{xx} {u}_{x}^{\varepsilon } + \frac{1}{2}{u}_{x}^{2}   {u}_{x}^{\varepsilon } = ( {{u}^{2} + \frac{1}{2}{\rho }^{2} - \Pi}) {u}_{x}^{\varepsilon }, \\  u^\varepsilon _{tx} {u}_{x} + {( u  {u}_{xx}) }^{\varepsilon}  {u}_{x} + \frac{1}{2}{( {u}_{x}^{2}) }^{\varepsilon}   {u}_{x} = [ {({u}^{2})^\varepsilon + \frac{1}{2}({\rho }^{2})^\varepsilon - \Pi^\varepsilon}] {u}_{x}.\end{array}\right.\]
Multiplying $(3.4)_3$ and $(3.5)_3$ by $\rho^\varepsilon$ and $\rho$ respectively, we have
\[\left\{  \begin{array}{l}   \rho_{t} \rho^\varepsilon+(\rho u)_{x}\rho^ \varepsilon=0 , \\ \rho_{t}^\varepsilon \rho+ (\rho u)_{x}^\varepsilon \rho=0 . \end{array}\right.\]
Summing the six equations above yields
\begin{equation}
\begin{split}
& ( u u^{\varepsilon} )_t + ( u_x u_x^{\varepsilon} )_t + ( \rho \rho^{\varepsilon} )_t
+ \underbrace{u^{\varepsilon} u u_x + u ( u u_x )^{\varepsilon} -u_x^{\varepsilon} u^2-u_x ( u^2 )^{\varepsilon}}_{\mathrm{I}} \\
&+ \underbrace{u_x^{\varepsilon} u u_{xx} + \frac{1}{2} u_x^{\varepsilon} u_x^2 + u_x ( u u_{xx} )^{\varepsilon} + \frac{1}{2} u_x ( u_x^2 )^{\varepsilon}}_{\mathrm{II}}
+u^{\varepsilon} \Pi_x +( u \Pi_x^{\varepsilon})\\
&+ \underbrace{\rho^{\varepsilon} ( \rho u )_x + \rho ( \rho u )_x^{\varepsilon}-\frac{1}{2} u_x^{\varepsilon} \rho^2 -u_x ( \frac{1}{2} \rho^2 )^{\varepsilon}}_{\mathrm{III}} \\
&=  - u_x^{\varepsilon} \Pi - u_x \Pi^{\varepsilon}.
\end{split}
\end{equation}
Rewriting equation $(3.6)$ in terms of $\mathrm{I}$, $\mathrm{II}$ and $\mathrm{III} $, we obtain\\
\[{\left( u{u}^{\varepsilon }\right) }_{t} + {\left( {u}_{x}{u}_{x}^{\varepsilon }\right) }_{t} + {\left( \rho {\rho}^{\varepsilon }\right) }_{t} + \mathrm{I} + \mathrm{II} + \mathrm{III} + {\left( {u}^{\varepsilon }\Pi \right)}_{x} + {\left( u {\Pi}^{\varepsilon }\right) }_{x} = 0.\]
To prove $(2.2)$, we need to regroup $(3.6)$, thanks to the Leibniz product rule, we have \\
$${u}^{\varepsilon }u{u}_{x} + u{\left( u{u}_{x}\right) }^{\varepsilon } = u^\varepsilon(\frac{1}{2}u^2)_x+u(\frac{1}{2}u^2)^\varepsilon_x= \frac{1}{2}{u}^{\varepsilon }({u}^{2})_{x} + \frac{1}{2}{u}({u^2})^\varepsilon_{x}$$
$$= \frac{1}{2}{( {u}^{\varepsilon } {u}^{2}) }_{x} - \frac{1}{2}{u}_{x}^{\varepsilon } {u}^{2} + \frac{1}{2}{[ u{( {u}^{2}) }^{\varepsilon }] }_{x} - \frac{1}{2}{u}_{x}{( {u}^{2}) }^{\varepsilon }.$$
Next, we proceed to calculate the terms $\mathrm{I}$, $\mathrm{II}$ and $\mathrm{III} $ individually,\\
\begin{equation}
\begin{split}
\mathrm{I} &=u^{\varepsilon} u u_x + u ( u u_x)^{\varepsilon} -u_x^{\varepsilon} u^2-u_x (u^2)^{\varepsilon}\\
&= \frac{1}{2} ( u^{\varepsilon} u^2 )_x
- \frac{1}{2} u_x^{\varepsilon} u^2
+ \frac{1}{2} [ u ( u^2 )^{\varepsilon} ]_x
- \frac{1}{2} u_x ( u^2 )^{\varepsilon}
- u_x^{\varepsilon} u^2
- u_x ( u^2 )^{\varepsilon} \\
&= \frac{1}{2} ( u^{\varepsilon} u^2 )_x
+ \frac{1}{2} [ u ( u^2 )^{\varepsilon} ]_x
- \frac{3}{2} u_x^{\varepsilon} u^2
- \frac{3}{2} u_x ( u^2 )^{\varepsilon},
\end{split}
\end{equation}
\begin{equation}
\begin{split}
\mathrm{II} &= u_x^{\varepsilon} u u_{xx}
+ \frac{1}{2} u_x^{\varepsilon}  u_x^2
+ u_x  ( u u_{xx} )_x^{\varepsilon}
+ \frac{1}{2} u_x ( u_x^2 )^{\varepsilon} \\
   &= u_x^{\varepsilon} (u u_x)_x-u_x^{\varepsilon} u^2_x
+ \frac{1}{2} u_x^{\varepsilon}  u_x^2
+ u_x  ( u u_x )_x^{\varepsilon}-u_x (u^2_x)^\varepsilon
+ \frac{1}{2} u_x  ( u_x^2 )^{\varepsilon} \\
&=u_x^{\varepsilon} (u u_x)_x-\frac{1}{2} u_x^{\varepsilon}  u_x^2
+ u_x  ( u u_x )_x^{\varepsilon}- \frac{1}{2} u_x ( u_x^2 )^{\varepsilon}\\
   &= [ u_x^{\varepsilon}  ( u u_x ) ]_x
- u_{xx}^{\varepsilon} ( u u_x )
+ [ u_x  ( u u_x )^{\varepsilon} ]_x \\
   &\quad - u_{xx}  ( u u_x )^{\varepsilon}
- \frac{1}{2} u_x^{\varepsilon}  u_x^2
- \frac{1}{2} u_x  ( u_x^2)^{\varepsilon},
\end{split}
\end{equation}
\begin{equation}
\begin{split}
\mathrm{III} &= \rho^{\varepsilon} ( \rho u )_x
+ \rho ( u \rho )_x^{\varepsilon}
- \frac{1}{2} u_x^{\varepsilon} \rho^2
- u_x  ( \frac{1}{2} \rho^2 )^{\varepsilon} \\
&= [ \rho^{\varepsilon} ( \rho u ) ]_x
- \rho_x^{\varepsilon}  ( \rho u )
+ [ \rho ( u \rho )^{\varepsilon} ]_x
- \rho_x ( u \rho )^{\varepsilon}
- \frac{1}{2} u_x^{\varepsilon}  \rho^2
- \frac{1}{2} u_x ( \rho^2 )^{\varepsilon}.
\end{split}
\end{equation}
Substituting $(3.7)-(3.9)$ into $(3.6)$, we obtain
\[{( u{u}^{\varepsilon}) }_{t} + {( {u}_{x}{u}_{x}^{\varepsilon}) }_{t} + {( \rho{\rho}^{\varepsilon}) }_{t} + \frac{1}{2}{( {u}^{\varepsilon }{u}^{2}) }_{x} + \frac{1}{2}{[ u{( {u}^{2}) }^{\varepsilon }] }_{x} - \frac{3}{2}{u}_{x}^{\varepsilon }   {u}^{2}  - \frac{3}{2}{u}_{x}  ({u}^{2})^{\varepsilon }\]
\[+ {( {u}_{x}^{\varepsilon }   u{u}_{x}) }_{x} - {u}_{xx}^{\varepsilon }u{u}_{x} + {[ {u}_{x}   {( u{u}_{x}) }^{\varepsilon }] }_{x} - {u}_{xx}{( u{u}_{x}) }^{\varepsilon } - \frac{1}{2}{u}_{x}^{\varepsilon }   {u}_{x}^{2} - \frac{1}{2}{u}_{x}   {( {u}_{x}^{2}) }^{\varepsilon }- {\rho }_{x}^{\varepsilon }   {\rho u}\]
\[+ {[ {\rho }^{\varepsilon }   ( \rho u) ] }_{x} + {[ \rho {( u\rho ) }^{\varepsilon }] }_{x} - {\rho }_{x}{( u\rho ) }^{\varepsilon } - \frac{1}{2}{u}_{x}^{\varepsilon }   {\rho }^{2} - \frac{1}{2}{u}_{x}  {( {\rho }^{2}) }^{\varepsilon } +(u^\varepsilon \Pi)_{x}+(u \Pi^\varepsilon)_{x}=0,\]
it follows that\\
\begin{equation}
\begin{split}
&(u u^\varepsilon)_t + (u_x u_x^\varepsilon)_t + (\rho \rho^\varepsilon)_t
+ \frac{1}{2} (u^\varepsilon u^2)_x + \frac{1}{2} [u(u^2)^\varepsilon]_x
+ (u_x^\varepsilon u u_x)_x  \\
&+[u_x(u u_x)^\varepsilon]_x+ [\rho^\varepsilon (\rho u)]_x + [\rho(u\rho)^\varepsilon]_x
+ (u^\varepsilon \Pi)_x + (u \Pi^\varepsilon)_x \\
=& \ \frac{3}{2} u_x^\varepsilon u^2 + \frac{3}{2} u_x (u^2)^\varepsilon
+ u_{xx}^\varepsilon u u_x + u_{xx}(u u_x)^\varepsilon
+ \frac{1}{2} u_x (u_x^2)^\varepsilon + \frac{1}{2} u_x^\varepsilon u_x^2 \\
&+ \rho_x^\varepsilon \rho u + \rho_x (u\rho)^\varepsilon
+ \frac{1}{2} u_x^\varepsilon \rho^2 + \frac{1}{2} u_x (\rho^2)^\varepsilon.
\end{split}
\end{equation}
For the sake of simplicity, we define $\delta u(l)=u(x+l)-u(x)=f-u$, $\delta u_x(l)=u_x(x+l)-u_x(x)=f _x-u_ x $ and $\delta \rho(l)=\rho(x+l)-\rho(x)=g-\rho$, by use of them we can abtain\\
\begin{equation}
\begin{split}
\int_{\mathbb{S}} \varphi_{\varepsilon}'(l) [ \delta u(l) ]^{3} dl
&= \int_{\mathbb{S}} \varphi_{\varepsilon}'(l) [ u(x+l) - u(x) ]^{3} dl \\
&= \int_{\mathbb{S}} \varphi_{\varepsilon}'(l) [ u^{3}(x+l) - 3u^{2}(x+l)u(x) \\
&\quad + 3u(x+l)u^{2}(x) - u^{3}(x) ] dl.
\end{split}
\end{equation}
By performing the change of variables $l=\gamma-x$ and carrying out direct computation, we have\\
\[
\begin{aligned}
\int_{\mathbb{S}} \varphi_{\varepsilon}'(l) u^{3}(x+l) dl
&= \int_{\mathbb{S}} \partial_{\gamma} \varphi_{\varepsilon}(\gamma - x) u^{3}(\gamma) d\gamma \\
&= - \int_{\mathbb{S}} \partial_{x} \varphi_{\varepsilon}(\gamma - x) u^{3}(\gamma) d\gamma \\
&= - \partial_{x} \int_{\mathbb{S}} \varphi_{\varepsilon}(\gamma - x) u^{3}(\gamma) d\gamma \\
&= - ( \varphi_{\varepsilon} * u^{3} )_{x}
= - ( u^{3} )_{x}^{\varepsilon}.
\end{aligned}
\]
Similarly, we get\\
\[
\begin{aligned}
\int_{\mathbb{S}} \varphi_{\varepsilon}'(l) [ -3u^{2}(x+l)u(x) ] dl
&= -3u(x) \int_{\mathbb{S}} \varphi_{\varepsilon}'(l) u^{2}(x+l) dl \\
&= [ -3u(x) ] [ - ( u^{2} )_{x}^{\varepsilon} ]
= 3u(x) ( u^{2} )_{x}^{\varepsilon},
\end{aligned}
\]

\[
\begin{aligned}
\int_{\mathbb{S}} \varphi_{\varepsilon}'(l) \left[ 3u(x+l)u^{2}(x) \right] dl
&= 3u^{2}(x) \int_{\mathbb{S}} \varphi_{\varepsilon}'(l) u(x+l) dl \\
&= -3u^{2} u^{\varepsilon}_x,
\end{aligned}
\]
\[\int_{\mathbb{S}}{\varphi }_{\varepsilon }'( l) u^3(x) {dl}=0.\]
Substituting these four identities into $(3.11)$, we obtain\\
\begin{equation}
\begin{aligned}
\int_{\mathbb{S}} \varphi_{\varepsilon}'(l) [ \delta u(l) ]^{3} dl
&= - ( u^{3} )_{x}^{\varepsilon}
+ 3u ( u^{2} )_{x}^{\varepsilon}
- 3u^{2} u_{x}^{\varepsilon} \\
&= - ( u^{3} )_{x}^{\varepsilon}
+ 3u ( u^{2} )_{x}^{\varepsilon}
- 3u^{2} u_{x}^{\varepsilon} \\
&= - ( u^{3} )_{x}^{\varepsilon}
+ 3 [u(u^2)^{\varepsilon}]_x
- 3u_{x} ( u^{2} )^{\varepsilon}
- 3u^{2} u_{x}^{\varepsilon},
\end{aligned}
\end{equation}
thus\\
\begin{equation}
\frac{3}{2}u^\varepsilon_x u^2+\frac{3}{2}u_x (u^2)^\varepsilon=-\frac{1}{2} \int _{\mathbb{S}}{\varphi }_{\varepsilon}'(l)[ {\delta u}( {l}) ]^3{dl}-\frac{1}{2}(u^3)^\varepsilon_x+\frac{3}{2}[u({u^2})^\varepsilon)]_x.
\end{equation}
By suitable modification of the deduction of the above processes, we can get\\
\begin{equation}
\begin{aligned}
&\int_{\mathbb{S}} \varphi_{\varepsilon}'(l) \delta u(l)[ \delta u_x(l)]^{2} dl\\
=&\int_{\mathbb{S}} \varphi_{\varepsilon}'(l) [u(x+l)u_x^2(x+l)
- 2u(x+l)u_x(x+l)u_x(x)
+ u(x+l)u_x^2(x) \\
&
- u(x)u_x^2(x+l)
+ 2u(x)u_x(x+l)u_x(x)
- u(x)u_x^2(x)
] dl \\
=& - ( u u_x^2)_{x}^{\varepsilon}+ 2u_x (u u_x )_{x}^{\varepsilon}
- u_x^2 u_{x}^{\varepsilon}\\
&+ u (u_x^2 )_{x}^{\varepsilon}- 2u u_x (u_x )_{x}^{\varepsilon}
+ \int_{\mathbb{S}} \varphi_{\varepsilon}'(l) u u_x^2 dl \\
=& - (u u_x^2 )_{x}^{\varepsilon}+ 2u_x (u u_x )_{x}^{\varepsilon}- u_x^2 u_{x}^{\varepsilon}
+ u(u_x^2)_{x}^{\varepsilon}- 2u u_x (u_x )_{x}^{\varepsilon}.
\end{aligned}
\end{equation}
Since
\[2{u}_{x}{( u{u}_{x}) }_{x}^{\varepsilon }+{u} {( {u}_{x}^{2}) }_{x}^{\varepsilon }=2[(u u_x)^\varepsilon u_x]_x-2(u u_x)^\varepsilon  u_{xx}+[u({u_x^2})^{\varepsilon }]_x-u_x(u_x^2)^\varepsilon,\]
we obtain
\begin{equation}
\begin{aligned}
& \int_{\mathbb{S}} \varphi_{\varepsilon}'(l) [ u(x+l) - u(x) ] [ u_x(x+l) - u_x(x) ]^{2} dl \\
=& - ( u u_x^2 )_{x}^{\varepsilon}
+ 2 [ ( u u_x )^{\varepsilon} u_x ]_{x}
- 2 ( u u_x )^{\varepsilon} u_{xx}
+ [ u ( u_x^2)^{\varepsilon} ]_{x} \\
&- u_x ( u_x^2 )^{\varepsilon}
- u_x^2 u_{x}^{\varepsilon}
- 2 u u_x ( u_x )_{x}^{\varepsilon}.
\end{aligned}
\end{equation}
Rearranging terms, we have\\
\[
\begin{aligned}
& 2( u u_x )^{\varepsilon} u_{xx}
+ u_x ( u_x^2 )^{\varepsilon}
+ u_x^2 u_{x}^{\varepsilon}
+ 2u u_x ( u_x )_{x}^{\varepsilon} \\
=& - \int_{\mathbb{S}} \varphi_{\varepsilon}'(l) ( u(x+l) - u(x) ) ( u_x(x+l) - u_x(x) )^{2} dl \\
&- ( u u_x^2 )_{x}^{\varepsilon}
+ 2[ ( u u_x )^{\varepsilon} u_x ]_{x}
+ [ u( u_x^2 )^{\varepsilon} ]_{x},
\end{aligned}
\]
it follows that\\
\begin{equation}
\begin{aligned}
& ( u u_x )^{\varepsilon} u_{xx}
+ u u_x  (u_{x})_x ^{\varepsilon}
+ \frac{1}{2}  u_x ( u_x^2 )^{\varepsilon} + \frac{1}{2}u_x^2 u_{x}^{\varepsilon}  \\
=& -\frac{1}{2} \int_{\mathbb{S}} \varphi_{\varepsilon}'(l) \delta u(l) [ \delta u_x(l)]^{2} dl
- \frac{1}{2} ( u u_x^2 )_{x}^{\varepsilon}
+ [( u u_x )^{\varepsilon} u_x ]_{x}
+ \frac{1}{2} [ u ( u_x^2 )^{\varepsilon} ]_{x}.
\end{aligned}
\end{equation}
Repeating the above argument, we obtain\\
\begin{equation}
\begin{aligned}
& \int_{\mathbb{S}} \varphi_{\varepsilon}'(l) \delta u(l) [ \delta \rho(l) ]^{2} dl \\
=& \int_{\mathbb{S}} \varphi_{\varepsilon}'(l) ( f - u ) ( g - \rho )^{2} dl \\
=& \int_{\mathbb{S}} \varphi_{\varepsilon}'(l) ( f - u ) ( g^{2} - 2g\rho + \rho^{2} ) dl \\
=& \int_{\mathbb{S}} \varphi_{\varepsilon}'(l) ( f g^{2} - 2fg\rho + f\rho^{2} - ug^{2} + 2g\rho u - u\rho^{2} ) dl \\
=& - ( u\rho^{2} )_{x}^{\varepsilon}
+ 2\rho ( u\rho )_{x}^{\varepsilon}
- \rho^{2} u_{x}^{\varepsilon}
+ u ( \rho^{2} )_{x}^{\varepsilon}
- 2u\rho \rho_{x}^{\varepsilon}.
\end{aligned}
\end{equation}
A routine computation yields\\
\[{2\rho }{( u\rho ) }_{x}^{\varepsilon} + u{( {\rho }^{2}) }_{x}^{\varepsilon }
= 2[\rho(\rho u)^{\varepsilon}]_x - 2{\rho }_{x}  {( u{\rho }) }^{\varepsilon} - {u}_{x}{( {\rho}^{2}) }^\varepsilon+[u(\rho^2)^{\varepsilon}]_{x}.\]
Substituting this expression into $(3.17)$, we obtain\\
\[
\begin{aligned}
& \int_{\mathbb{S}} \varphi_{\varepsilon}'(l) \delta u(l) [ \delta \rho(l) ]^{2} dl \\
=& - ( u\rho^{2} )_{x}^{\varepsilon} + 2 [ \rho( u\rho )^{\varepsilon}]_{x} - 2\rho_{x} ( u\rho )^{\varepsilon} + [ u ( \rho^{2} )^{\varepsilon} ]_{x} - u_{x} ( \rho^{2} )^{\varepsilon} - \rho^{2} u_{x}^{\varepsilon} - 2u\rho \rho_{x}^{\varepsilon}.
\end{aligned}
\]
So\\
\begin{equation}
\begin{aligned}
& \rho_{x}^{\varepsilon} \rho u + \rho_{x} ( u\rho )^{\varepsilon} + \frac{1}{2} u_{x}^{\varepsilon} \rho^{2} + \frac{1}{2} u_{x} ( \rho^{2} )^{\varepsilon} \\
= &-\frac{1}{2} \int_{\mathbb{S}} \varphi_{\varepsilon}'(l) \delta u(l) [ \delta \rho(l) ]^{2} dl
- \frac{1}{2} ( u\rho^{2} )_{x}^{\varepsilon}
+ [ \rho ( u\rho )^{\varepsilon} ]_{x}
+ \frac{1}{2} [ u ( \rho^{2} )^{\varepsilon}]_{x}.
\end{aligned}
\end{equation}
Substituting $(3.13)$, $(3.16)$ and $(3.18)$ into $(3.10)$, we have\\
\[
\begin{aligned}
& (u u^\varepsilon)_t+(u_x u^\varepsilon_x)_t+(\rho \rho^\varepsilon)_t
+\frac{1}{2}(u^\varepsilon u^2)_x + \frac{1}{2}[u (u^2)^\varepsilon]_x+(u^\varepsilon_x u u_x)_x \\
& +[u_x (u u_x)^\varepsilon]_x+[\rho^\varepsilon(\rho u)]_x+[\rho(u \rho)^\varepsilon]_x +(u^\varepsilon \Pi)_x+(u \Pi^\varepsilon)_x \\
=& -\frac{1}{2}\int_{\mathbb{S}}\varphi_\varepsilon'(l)\bigl[\delta u(l)\bigr]^3 dl
-\frac{1}{2}(u^3)^\varepsilon_x+\frac{3}{2}[u(u^2)^\varepsilon]_x \\
&-\frac{1}{2}\int_{\mathbb{S}}\varphi_\varepsilon'(l)\delta u(l)[\delta u_x(l)]^2 dl
-\frac{1}{2}(u u_x^2)^\varepsilon_x+[(u u_x)^\varepsilon u_x]_x
+\frac{1}{2}[u(u_x^2)^\varepsilon]_x \\
&-\frac{1}{2}\int_{\mathbb{S}}\varphi_\varepsilon'(l)\delta u(l)[\delta \rho(l)]^2 dl
-\frac{1}{2}(u\rho^2)^\varepsilon_x+[\rho(u\rho)^\varepsilon]_x
+\frac{1}{2}[u(\rho^2)^\varepsilon]_x,
\end{aligned}
\]
it follows that\\
\begin{equation}
\begin{aligned}
& (u u^{\varepsilon})_{t} + (u_{x} u_{x}^{\varepsilon})_{t} + (\rho \rho^{\varepsilon})_{t}
+ \frac{1}{2} (u^{\varepsilon} u^{2} )_{x}
- [u( u^{2} )^{\varepsilon} ]_{x}
+ (u_{x}^{\varepsilon} u u_{x} )_{x} \\
&+[\rho^{\varepsilon} (\rho u) ]_{x}
+ (u^{\varepsilon} \Pi )_{x}
+ (u \Pi^{\varepsilon} )_{x}
+ \frac{1}{2} ( u^{3} )_{x}^{\varepsilon}
+ \frac{1}{2} ( u u_{x}^{2} )_{x}^{\varepsilon}
- \frac{1}{2} [ u ( u_{x}^{2} )^{\varepsilon} ]_{x} \\
&+ \frac{1}{2} ( u \rho^{2} )_{x}^{\varepsilon}
- \frac{1}{2} [ u ( \rho^{2} )^{\varepsilon} ]_{x} \\
=& -\frac{1}{2} \int_{\mathbb{S}} \varphi_{\varepsilon}'(l) [ \delta u(l) ]^{3} dl
- \frac{1}{2} \int_{\mathbb{S}} \varphi_{\varepsilon}'(l) \delta u(l) [ \delta u_{x}(l) ]^{2} dl \\
&\quad - \frac{1}{2} \int_{\mathbb{S}} \varphi_{\varepsilon}'(l) \delta u(l) [ \delta \rho(l) ]^{2} dl \\
= &D^{2CH}_{\varepsilon}(u, \rho).
\end{aligned}
\end{equation}
The proof of the first conclusion in Theorem $2.1$ is now complete.\ Next, we prove the second conclusion.\ To prove that $u u_x u^\varepsilon_x\rightarrow u u^2_x$ in the sense of distribution, we first verify that the nonlinear term $u u^2_x$ belongs to $L^{1}([0,T];L^{1}(\mathbb{S}))$, so that the following distributional calculation are well defined.\ For $\forall t \in R^{+}$ fixed, thanks to H\"{o}lder inequality, we infer that
\[
\begin{aligned}
& \int_{\mathbb{S}} | u_{x}^{2}(t,x) | | u(t,x) | dx
\leq \| u_{x}^{2} \|_{L^{\frac{2+\eta}{2}}(\mathbb{S})}\| u \|_{L^{\frac{2+\eta}{\eta}}(\mathbb{S})} \\[4pt]
=&( \int_{\mathbb{S}} u_{x}^{2+\eta} dx )^{\frac{2}{2+\eta}} ( \int_{\mathbb{S}} u^{\frac{2+\eta}{\eta}} dx )^{\frac{\eta}{2+\eta}}
= \| u_{x} \|_{L^{2+\eta}(\mathbb{S})}^{2} \| u \|_{L^{\frac{2+\eta}{\eta}}(\mathbb{S})}.
\end{aligned}
\]
Integrating both sides of this inequality with respect to $t$ over $[0,T]$, we obtain\\
\[
\begin{aligned}
& \int_{0}^{T} \int_{\mathbb{S}} | u_{x}^{2}(t,x) | | u(t,x) | dxdt\\
\leq& \int_{0}^{T} \| u_{x} \|_{L^{2+\eta}(\mathbb{S})}^{2} \| u \|_{L^{\frac{2+\eta}{\eta}}(\mathbb{S})} dt \\[4pt]
\leq &\ \| \| u_{x} \|^{2}_{L^{2+\eta}(\mathbb{S})} \|_{L^{\frac{2+\delta}{2}}[0,T]}
\| \| u \|_{L^{\frac{2+\eta}{\eta}}(\mathbb{S})} \|_{L^{\frac{2+\delta}{\delta}}[0,T]} \\[4pt]
\leq & ( \int_{0}^{T} \| u_{x} \|_{L^{2+\eta}(\mathbb{S})}^{2+\delta} dt )^{\frac{2}{2+\delta}}
\| \| u \|_{L^{\frac{2+\eta}{\eta}}(\mathbb{S})} \|_{L^{\frac{2+\delta}{\delta}}[0,T]} \\[4pt]
=& \| u_{x} \|_{L^{2+\delta}([0,T];L^{2+\eta}(\mathbb{S}))}^{2}
\| u \|_{L^{\frac{2+\delta}{\delta}}([0,T];L^{\frac{2+\eta}{\eta}}(\mathbb{S}))}.
\end{aligned}
\]
Thus, we obtain\\
\[\|u^2_x u\|_{{L^1}([0,T];L^1(\mathbb{S})) } \leq \|u_x\|^2_{{L^{2+\delta }}([0,T];L^{{2+\eta }}(\mathbb{S}))} \|u\|_{{L^\frac{2+\delta }{\delta  }}([0,T];L^{\frac{2+\eta }{\eta  }}(\mathbb{S}))}.\]
Next, we prove that $(u u_x)u^\varepsilon_x\rightarrow u u^2_x$, in the sense of distribution, that is for any test function $\psi \in C^\infty_0 ([0,T]\times \mathbb{S})$, one has \\
\[\mathop{\lim }\limits_{{\varepsilon  \rightarrow  0}}\int _{0}^{T}\int _{\mathbb{S}}(u u_x u^\varepsilon_x \psi-u u^2_x \psi)  {dxdt} =0.\]
By H\"{o}lder inequality, we have\\
\[
\begin{aligned}
& | \int_{0}^{T} \int_{\mathbb{S}} u u_{x} ( u_{x}^{\varepsilon} - u_{x} ) \psi \, dxdt |\\
 \leq& \int_{0}^{T} \int_{\mathbb{S}} |u u_{x} ( u_{x}^{\varepsilon} - u_{x} ) \psi |\, dxdt \\
\leq & \| \psi \|_{L^{\infty}([0,T]\times\mathbb{S})}
\int_{0}^{T} \int_{\mathbb{S}} | u u_{x} ( u_{x}^{\varepsilon} - u_{x} ) | dxdt \\
 \leq & \| \psi \|_{L^{\infty}([0,T]\times\mathbb{S})}
\| u_{x} \|_{L^{2+\delta}([0,T];L^{2+\eta}(\mathbb{S}))}
\| u_{x}^{\varepsilon} - u_{x} \|_{L^{2+\delta}([0,T];L^{2+\eta}(\mathbb{S}))} \\
&   \| u \|_{L^{\frac{2+\delta}{\delta}}([0,T];L^{\frac{2+\eta}{\eta}}(\mathbb{S}))}.
\end{aligned}
\]
In view of $u_x\in L^{2+}([0,T];L^{2+}(\mathbb{S}))$ and
\[\|u^\varepsilon_x-u_x\|_{{L^{2+\delta}}([0,T];L^{2+\eta}(\mathbb{S}))} \rightarrow 0,\]
we get $(u u_x)u^\varepsilon_x \rightarrow (u u_x)u_x$ in the sense of distributions.\

Additionally, the celebrated H\"{o}rmander-Mihlin multiplier theorem on $\mathbb{R}$ tells that if $u$, $u_x$ and $ \rho \in L^{2+}([0,T];L^{2+}(\mathbb{R}))$ and $-\Pi_{xx}+\Pi=u^2+\frac{1}{2}u^2_x+\frac{1}{2}\rho^2,$ then $\Pi\in L^{1+}([0,T];W^{2,1+}(\mathbb{R}))$.\ Lemma 2.2 can support  multiplier theorem transplanting from $\mathbb{R}$ to $\mathbb{S}$, guaranteeing when $u, {u}_{x} \in  {L}^{2 + }( {\lbrack  {0,T}\rbrack  ;{L}^{2 + }( \mathbb{S}) })$ in the text, $\Pi$ and other non-local terms can be approximated in the periodic domain and limits can be taken.\ Therefore, $\Pi\in L^{1+}([0,T];W^{2,1+}(\mathbb{S}))$.\ Then, we will prove that $\Pi u^\varepsilon$ and $\Pi^\varepsilon u$ converge to $\Pi u$ in the sense of distributions.\\
\[
\begin{aligned}
& | \int_{0}^{T} \int_{\mathbb{S}} ( \Pi^{\varepsilon} u - \Pi u ) \psi \, dxdt | \,  \\
\leq & \int_{0}^{T} \int_{\mathbb{S}} |( \Pi^{\varepsilon} u - \Pi u ) \psi| \, dxdt \\
\leq & \| \psi \|_{L^{\infty}([0,T] \times \mathbb{S})}
\int_{0}^{T} \int_{\mathbb{S}} ( \Pi^{\varepsilon} - \Pi ) u \, dxdt \\
\leq & \| \psi \|_{L^{\infty}([0,T] \times \mathbb{S})}
\| \Pi^{\varepsilon} - \Pi \|_{L^{\frac{2+\delta}{1+\delta}}(0,T;L^{\frac{2+\eta}{1+\eta}}(\mathbb{S}))}
\| u \|_{L^{2+\delta}([0,T];L^{2+\eta}(\mathbb{S}))}.
\end{aligned}
\]
Since $\Pi^\varepsilon\rightarrow \Pi$ in $L^{1+}([0,T];W^{2,1+}(\mathbb{S}))$, $\Pi^\varepsilon u\rightarrow \Pi u$ in the sense of distributions.
Subsequently, $u\in{{L^\infty}([0,T];H^1(\mathbb{S}))}$ and $u_x \in{{L^{2+}}([0,T];L^{2+}(\mathbb{S}))}$ imply that $(u^2)^\varepsilon u\rightarrow u^3$ and $u_x u^\varepsilon_x \rightarrow u^2_x$ in the sense of distributions.\ Finally, we obtain that $D^{2CH}_\varepsilon$ converges to
\[(u^2)_t + (u^2_x)_t+ (\rho^2)_t + (u u^2_x)_x+{( {\rho}^{2}u) }_{x} + {2}{( {\Pi}{u}) }_{x}\]
in the sense of distributions as $\varepsilon\rightarrow 0,$ which completes the proof of Theorem 2.1.

\end{proofof}

\section{Proofs of Corollaries 2.1 and Theorem 2.2}
\newtheorem{theorem4}{Theorem}[section]
\newtheorem{lemma4}{Lemma}[section]
\newtheorem {remark4}{Remark}[section]
\newtheorem{corollary4}{Corollary}[section]

As an application of the convergence analysis for the regularized dissipative term $D^{2CH}_\varepsilon(u,\rho)$ established in Section 3, we derive a corollary and a theorem regarding the structural properties of weak solutions for system $(2.1)$.\ These results establish sufficient conditions that guarantee the energy conservation property of weak solutions for system $(2.1)$.\ Especially, Theorem $2.2$ reveals that the Onsager exponent of system $(2.1)$ is 1, with the solution being unique.\ Firstly, we are in the position to prove Corollary $2.1$.
\vspace{1em}
\begin{proofof}{Corollary 2.1}
According to Fubini$'$s theorem, one has\\
\begin{equation}
\int _{\mathbb{S}}\left| {\int _{\mathbb{S}}{\varphi }_{\varepsilon }'( l) { \delta u( l )  }^{3}dl }\right| {dx} \leq  \int _{\mathbb{S}}| {{\varphi }_{\varepsilon }'( l ) }| \left( {\int _{\mathbb{S}}{| \delta u( l ) | }^{3}{dx}}\right) dl.
\end{equation}
By $(2.4)$, we have\\
\begin{equation}
\int _{\mathbb{S}}\left| {\int _{\mathbb{S}}{\varphi }_{\varepsilon }'( l ) { \delta u( l)  }^{3}dl }\right| {dx} \leq  \int _{\mathbb{S}}\left| {{\varphi }_{\varepsilon }'( l ) }\right| {C}_{1}( t) | l | {\sigma }_{1}( l ) dl.
\end{equation}
Integrating both sides of (4.2) with respect to t and letting $\frac{l}{\varepsilon}=\xi$, noticing that ${C}_{1}(t)$ is an integrable function on $[0,T]$ and $\sigma_1(l)$ is a bounded function on some neighborhood of the origin, we deduce that\\
\begin{equation}
\begin{aligned}
&\int_{0}^{T}\int_{\mathbb{S}}
\left| \int_{\mathbb{S}} \varphi_{\varepsilon}'(l) [ \delta u(l) ]^{3} dl \right| dxdt\\
\leq \ & \left( \int_{0}^{T} C_{1}(t) dt \right) \int_{\mathbb{S}} | \varphi_{\varepsilon}'(l) | |l| \sigma_{1}(l) dl  \\
\leq \ & C \left( \int_{|\xi|<1} | \varphi'(\xi)| |\xi| \sigma_{1}(\varepsilon \xi) d\xi \right)
\end{aligned}
\end{equation}
By an argument analogous to that leading to $(4.3)$, we obtain\\
\begin{equation}
\begin{aligned}
&\int_{0}^{T}\int_{\mathbb{S}}
\left| \int_{\mathbb{S}} \varphi_{\varepsilon}'(l)\delta u(l) [ \delta u_x(l) ]^{2} dl \right| dxdt\\
\leq \ & \left( \int_{0}^{T} C_{2}(t) dt \right) \int_{\mathbb{S}} | \varphi_{\varepsilon}'(l) | |l| \sigma_{2}(l) dl  \\
\leq \ & C \left( \int_{|\xi|<1} | \varphi'(\xi)| |\xi| \sigma_{2}(\varepsilon \xi) d\xi \right)
\end{aligned}
\end{equation}
and
\begin{equation}
\begin{aligned}
&\int_{0}^{T}\int_{\mathbb{S}}
\left| \int_{\mathbb{S}} \varphi_{\varepsilon}'(l)\delta u(l) \left[ \rho(l) \right]^{2} dl \right| dxdt \\
\leq \ & \left( \int_{0}^{T} C_{3}(t) dt \right) \int_{\mathbb{S}} \left| \varphi_{\varepsilon}'(l) \right| |l| \sigma_{3}(l) dl  \\
\leq \ & C \left( \int_{|\xi|<1} \left| \varphi'(\xi) \right| |\xi| \sigma_{3}(\varepsilon \xi) d\xi \right).
\end{aligned}
\end{equation}
Finally, recalling that $\sigma_i(l)\rightarrow 0$ as $l\rightarrow 0$, so $\sigma_i(\varepsilon \xi)\rightarrow 0$ as $\varepsilon\rightarrow 0$, and noting that $\sigma_i(l)$ is bounded on some neighbourhood of the origin while $|\varphi_{\varepsilon}'(\xi)||\xi|$ is integrable for $|\xi|<1$, the lebesgue dominated convergence theorem applies.\ Therefore, we may pass the limit inside the integral, yielding\\
\begin{equation}
\begin{aligned}
\mathop{\lim }\limits_{{\varepsilon  \rightarrow  0}}\int _{| \xi|  \leq   1  } |\varphi_{\varepsilon}'(l)||\xi| \sigma_i(\varepsilon \xi){d\xi}                   = \int _{| \xi|  \leq 1}(\mathop{\lim }\limits_{{\varepsilon  \rightarrow  0}}|\varphi_{\varepsilon}'(l)||\xi|\sigma_i(\varepsilon \xi) ){d\xi} = 0.
\end{aligned}
\end{equation}
By integrating $(2.2)$ both in time over $(0,t)$ and space over $\mathbb{S}$, we have
\begin{equation}
\int_{0}^{t}\int_{\mathbb{S}}(u u^\varepsilon)_t +(u_x u^\varepsilon_x)_t +(\rho \rho^\varepsilon)_t dxdt= \int _{0}^{t}\int _{\mathbb{S}}{D}_{\varepsilon }^{2CH}( {u,\rho }) {dxdt}.
\end{equation}
Since
\[\left|\int _{0}^{t}\int _{\mathbb{S}}{D}_{\varepsilon }^{2CH}( {u,\rho }) {dxdt}\right|  \leq  \int _{0}^{t}\int_{\mathbb{S}}|{D}_{\varepsilon }^{2CH}( {u,\rho }) |{dxdt},\]
it follows from $(4.3)-(4.6)$ that $\int _{0}^{t}\int_{\mathbb{S}}{D}_{\varepsilon }^{2CH}( {u,\rho }){dxdt} \rightarrow 0$ as $\varepsilon \rightarrow 0$.\
Letting $\varepsilon \rightarrow 0$ on both side of $(4.7)$, we have
\[\int_{0}^{t} \frac{d}{dt}\int_{\mathbb{S}}(u^2 + u_x^2 +\rho^2)dxdt=0,\]
that is
\[\int _{\mathbb{S}}( {{u}^{2} + {u}_{x}^{2} + {\rho }^{2}}) {dx} = \int _{\mathbb{S}}[ {{u}_{0}^{2} + (u_0)_{x}^{2} + {\rho }_0^{2}}] {dx}.\]
This completes the proof of Corollary $2.1$.
\end{proofof}
\vspace{1em}

\begin{proofof}{Theorem 2.2}\\
\\
{\bf Step 1} \ The case of $p=q=3$\\
\

According to H\"{o}lder inequality and the mean value theorem, we get\\
\begin{equation}
\left| {\int _{\mathbb{S}}{ \delta u( l )   }^{3}{dx}}\right|  \leq  \| {\delta u}(l )\|_{{L}^{3}( \mathbb{S}) }^{2}\| {\delta u}( l )\|_{{L}^{3}( \mathbb{S}) }.
\end{equation}
and\\
\begin{equation}
\|{\delta u}( l )\|_{{L}^{3}( \mathbb{S}) } \leq |l|\| u\| _{{W}^{1,3}( \mathbb{S}) },
\end{equation}
Inserting $(4.9)$ into $(4.8)$, we arrive at\\
\begin{equation}
\left| {\int _{\mathbb{S}}{  \delta u(l) }^{3}{dx}}\right| \leq  |l|\| {\delta u}(l )\|_{{L}^{3}( \mathbb{S}) }^{2} \| u{\| }_{{W}^{1,3}( \mathbb{S}) },
\end{equation}
Using Fubini$'$s theorem and $(4.10)$, we get\\
\begin{equation}
\begin{aligned}
\int _{\mathbb{S}}\left| {\int _{\mathbb{S}}{\varphi }_{\varepsilon }'( l ) { \delta u( l ) }^{3}dl }\right| {dx}
&\leq  \int _{\mathbb{S}}| {{\varphi }_{\varepsilon }'( l ) }| \left( {{\int }_{\mathbb{S}}{| \delta u( l ) | }^{3}{dx}}\right) dl\\
&\leq  \int _{\mathbb{S}}| {{\varphi }_{\varepsilon }'( l ) }| | l | \| {\delta u}{\| }_{{L}^{3}( \mathbb{S}) }^{2}d{l \| u{\| }_{{W}^{1,3}( \mathbb{S}) }}.
\end{aligned}
\end{equation}
By means of the change of variables $\eta=\frac{l}{\varepsilon}$, we can reformulate (4.11) as\\
\begin{equation}
\begin{aligned}
&\int _{\mathbb{S}}\left| {\int _{\mathbb{S}}{\varphi }_{\varepsilon }'( l) { \delta u( l)   }^{3}{dl}}\right| {dx} \\
\leq \ & \int _{| \eta |  \leq  1}| {{\varphi }'( \eta ) }| | \eta | \| u( {x + {\varepsilon \eta },t})  - u( {x,t}) {\| }^{2}_{{L}^{3}(\mathbb{S}) }{d\eta }\| u{\| }_{{W}^{1,3}(\mathbb{S}) }.
\end{aligned}
\end{equation}
Next, integrate both sides of equation $(4.12)$ with respect to t, we obtain\\
\begin{equation}
\begin{aligned}
&\int _{0}^{T}\int _{\mathbb{S}}\left| {\int _{\mathbb{S}}{\varphi }_{\varepsilon }'( l ) {  \delta u( l ) }^{3}dl }\right| {dxdt} \\
\leq \ & \int _{| \eta |  \leq  1}| {{\varphi }'( \eta ) }| | \eta | \| u( {x + {\varepsilon \eta },t})  - u( {x,t}) {\| }_{{L}^{3}( {  [0,T] ;{L}^{3}( \mathbb{S}) }) }^{2}d\eta \| u{\| }_{{L}^{3}( { [0,T]  ;{W}^{1,3}( \mathbb{S}) }) }.
\end{aligned}
\end{equation}
For fixed \(\eta\), as \(\varepsilon \rightarrow0\), we have \(\varepsilon\eta\rightarrow0\). By strong continuity of translations in lebesgue spaces and the fact that $u\in L^3([0,T];L^3(\mathbb{S}))$, we have$\|u(x+\varepsilon\eta,t)-u(x,t)\|_{L^3([0,T];L^3(\mathbb{S}))}\to 0$. Since $\|u(x+\varepsilon\eta,t)-u(x,t)\|^2_{L^3([0,T];L^3(\mathbb{S}))}\le 4\|u\|^2_{L^3([0,T];L^3(\mathbb{S}))}$ and $|\varphi'(\eta)||\eta|$ is integrable over the domain of integration, the lebesgue dominated convergence theorem applies.\ Therefore, we can pass the limit inside the integral, yielding
\begin{equation}
\begin{aligned}
&\mathop{\lim }\limits_{{\varepsilon  \rightarrow  0}}\int _{|\eta|\leq1}|{\varphi }'( \eta )||\eta|||u(x+\varepsilon\eta,t)-u(x,t)||^2_{{L}^{3}( { [0,T] ;{L}^{3}( \mathbb{S}) }) } d\eta\\
&  \cdot ||u|| _{{L}^{3}( { [0,T] ;{W}^{1,3}( \mathbb{S}) }) }  \\
= \ & \int _{|\eta|\leq1}\mathop{\lim }\limits_{{\varepsilon  \rightarrow  0}} |{\varphi }'( \eta )||\eta|||u(x+\varepsilon\eta,t)-u(x,t)||^2_{{L}^{3}( {  [0,T] ;{L}^{3}( \mathbb{S}) }) } d\eta \\
& \cdot ||u|| _{{L}^{3}( { [0,T]  ;{W}^{1,3}( \mathbb{S}) }) } \\
= \ & 0,
\end{aligned}
\end{equation}
we find\\
\begin{equation}
\mathop{\lim }\limits_{{\varepsilon  \rightarrow  0}}\int _{\mathbb{S}}{\varphi }_{\varepsilon }'( l ) [{ \delta u( l )  }]^{3}dl  = 0.
\end{equation}
By an argument analogous to that leading to $(4.10)$, we get\\
\begin{equation}
\int _{\mathbb{S}}{\delta u}( l ) {[ \delta {u}_{x}( l ) ]  }^{2}{dx} \leq  | l | \|\delta {u}_{x}\|_{{L}^{3}( \mathbb{S}) }^{2}\| u\| _{{W}^{1,3}( \mathbb{S}) }.
\end{equation}
\begin{equation}
\int _{\mathbb{S}}{\delta u}( l ) {[  \delta \rho( l ) ]  }^{2}{dx} \leq  | l | \|\delta \rho(l)\|_{{L}^{3}( \mathbb{S}) }^{2}\| u{\| }_{{W}^{1,3}( \mathbb{S}) }.
\end{equation}
Similar to get $(4.13)$, we obtain\\
\begin{equation}
\begin{aligned}
& \int_{0}^{T} \int_{\mathbb{S}} \left| \int_{\mathbb{S}} \varphi_{\varepsilon}'(l) \delta u(l) [\delta u_{x}(l)]^{2} dl \right| dxdt  \\
\leq \ &  \int_{|\eta| \leq 1} |\varphi'(\eta)| |\eta| \| u_{x}(x+\varepsilon\eta,t) - u_{x}(x,t) \|_{L^{3}([0,T];L^{3}(\mathbb{S}))}^{2} d\eta\ \\
&  \ \ \ \cdot \| u \|_{L^{3}([0,T];W^{1,3}(\mathbb{S}))}
\end{aligned}
\end{equation}
and
\begin{equation}
\begin{aligned}
&\int _{0}^{T}\int _{\mathbb{S}}\left| {\int _{\mathbb{S}}{\varphi }_{\varepsilon }'( l ) {\delta u}( l ) {\lbrack  \delta \rho( l ) \rbrack  }^{2}dl }\right| {dxdt} \\
&\leq \int _{| \eta |  \leq  1}| {{\varphi }'( \eta ) }| | \eta | \| \rho( {x + {\varepsilon \eta },t})  - \rho( {x,t}) \|_{{L}^{3}( {[  {0,T}]  ;{L}^{3}( \mathbb{S}) }) }^{2}d\eta \\
& \ \ \ \ \  \cdot \| u{\| }_{{L}^{3}( {[  0,T] ;{W}^{1,3}( \mathbb{S}) }) }
\end{aligned}
\end{equation}
Using the strong continuity of the translation operator on the Lebesgue spaces, $u \in  {L}^{3}\left( {[0,T];{W}^{1,3}\left( \mathbb{S}\right) }\right)$ and  dominated convergence theorem, it follows that\\
\begin{equation}
\mathop{\lim }\limits_{{\varepsilon  \rightarrow  0}}\int _{\mathbb{S}}{\varphi }_{\varepsilon }'( l ) \delta u(l) {[  \delta u_x( l ) ] }^{2}dl  = 0.
\end{equation}
\begin{equation}
\mathop{\lim }\limits_{{\varepsilon  \rightarrow  0}}\int _{\mathbb{S}}{\varphi }_{\varepsilon }'( l ) \delta u(l) {[ \delta \rho( l ) ]  }^{2}dl  = 0.
\end{equation}
By integrating $(2.2)$ both in time over $(0,t)$ and space over $\mathbb{S}$, using the characteristic of periodic function and letting $\varepsilon\rightarrow 0$, we obtain
$$\int _{\mathbb{S}}( {{u}^{2} + {u}_{x}^{2} + {\rho }^{2}}) {dx} = \int _{\mathbb{S}}[ {{u}_{0}^{2} + (u_0)_{x}^{2} + {\rho }_0^{2}}] {dx}.$$
This means that energy is conserved.
\\
\\
{\bf Step 2} The case of $\frac{1}{p}+\frac{2}{q}=1.$
\\

First, we prove that $(u_{x},\rho)\in L^{3}([0,T], L^{3}(\mathbb{S})\times L^{3}(\mathbb{S})).$ Note that,
\[\|(u_x,\rho)\|^3_{{L}^{3}\left( {\left\lbrack  {0,T}\right\rbrack  ;L^3 (\mathbb{S})\times{L}^{3}\left( \mathbb{S}\right) }\right)  }=\|u_x\|^3_{{L}^{3}\left( {\left\lbrack  {0,T}\right\rbrack  ;L^3 (\mathbb{S}) }\right)  }+\|\rho\|^3_{{L}^{3}\left( {\left\lbrack  {0,T}\right\rbrack  ;L^3 (\mathbb{S}) }\right)  },\]
we only need to prove $u_x\in{{L}^{3}( {[  {0,T}]  ;L^3 (\mathbb{S}) })  }$ and $\rho\in{{L}^{3}( {[  {0,T}]  ;L^3 (\mathbb{S}) })  }$.
By the interpolation inequality, for any measurable function $f$,
\[\|f\|_{{L}^{p}( \mathbb{S}) } \leq  C\|f\|_{{L}^{p_{0}}( \mathbb{S}) }^{1-\theta} \|f\|_{{L}^{p_1}( \mathbb{S}) }^{\theta},\]
where $1\leq p$, $p_0$, $p_1\leq\infty$, and $\theta\in(0,1)$ satisfies $\frac{1}{p}=\frac{1-\theta}{p_0}+\frac{\theta}{p_1}$, we have
\begin{equation}
\|u_x\|_{{L}^{3}( \mathbb{S}) } \leq  C\|u_x\|_{{L}^{2}( \mathbb{S}) }^{\frac{{2q} - 6}{{3q} - 6}}\|u_x\|_{{L}^{q}( \mathbb{S}) }^{\frac{q}{{3q} - 6}},
\end{equation}
then
\begin{equation}
\|u_x\|^3_{{L}^{3}(\mathbb{S}) } \leq  C \|u_x\|_{{L}^{2 }(\mathbb{S})}^{\frac{{2q} - 6}{q - 2}}\|u_x\|_{{L}^{q}( \mathbb{S}) }^{\frac{q}{q - 2}},
\end{equation}
it follows that
\begin{equation}
\begin{split}
\| u_x \|_{L^3([0,T];L^3(\mathbb{S}))}^3
&= \left( \int_0^T \| u_x \|_{L^3(\mathbb{S})}^3 dt \right)^{\frac{1}{3}} \\
&\leq C \left( \int_0^T \| u_x \|_{L^2(\mathbb{S})}^{\frac{2q-6}{q-2}} \| u_x \|_{L^q(\mathbb{S})}^{\frac{q}{q-2}} dt \right)^{\frac{1}{3}} \\
&\leq C \| u_x \|_{L^\infty([0,T];L^2(\mathbb{S}))}^{\frac{2q-6}{q-2}}
\| u_x \|_{L^{\frac{q}{q-2}}([0,T];L^q(\mathbb{S}))}^{\frac{q}{3(q-2)}}.
\end{split}
\end{equation}
Since $u_x\in {L^{p}([0,T];L^q(\mathbb{S}))}$ with $\frac{1}{p}+\frac{2}{q}=1$,
we have $u_x\in{L^{3}([0,T];L^3(\mathbb{S}))}$.

Similarly, following the derivation of $(4.24)$, we arrive at
\begin{equation}
\begin{split}
\| \rho \|_{L^3([0,T];L^3(\mathbb{S}))}^3
&= \left( \int_0^T \| \rho \|_{L^3(\mathbb{S})}^3 dt \right)^{\frac{1}{3}} \\
&\leq C \left( \int_0^T \| \rho \|_{L^2(\mathbb{S})}^{\frac{2q-6}{q-2}} \| \rho \|_{L^q(\mathbb{S})}^{\frac{q}{q-2}} dt \right)^{\frac{1}{3}} \\
&\leq C \| \rho \|_{L^\infty([0,T];L^2(\mathbb{S}))}^{\frac{2q-6}{q-2}}
\| \rho \|_{L^{\frac{q}{q-2}}([0,T];L^q(\mathbb{S}))}^{\frac{q}{3(q-2)}}.
\end{split}
\end{equation}
Since $\rho\in {L^{p}([0,T];L^q(\mathbb{S}))}$ satisfies $\frac{1}{p}+\frac{2}{q}=1$, we have
$\rho\in{L^{3}([0,T];L^3(\mathbb{S}))}$. It follows from Step 1 that the energy of the weak solution $(u,\rho)$  is conserved.

Next, we prove the uniqueness of the weak solution.
In fact, this process is very similar to the proof of Theorem 2.1. For the sake of completeness, we include it here.

In view of the regularity of the solutions, we need to perform the operations in the sense of distributions. For
$\psi\in C_{0}^{\infty}([0,T]\times\mathbb{S}),$ by $(3.2)_1$, we have
\[\langle u_{t}+uu_{x}+\Pi_{x}, \psi_{x}\rangle=0,\]
it follows that
\[\langle u_{tx}+u_{x}^{2}+u u_{xx}+\Pi_{xx}, \psi\rangle=0,\]
that is
\[u_{tx}+u_{x}^{2}+u u_{xx}+\Pi_{xx}=0\] in the sense of distributions. The subsequent calculations are carried out in the sense of distributions.
For convenience, we omit the explicit writing of the test functions.

Substituting the expression for $\Pi_{xx}$ given in $(3.3)$ into this identity and combining it with $(3.2)_1$ and $(3.2)_3$, we obtain \\
\begin{equation}
\left\{\begin{array}{ll}
{u}_{t} + u{u}_{x} + {\Pi}_{x} = 0,\\
 {u}_{tx} + u{u}_{xx} + \frac{1}{2}{u}_{x}^{2} = {u}^{2} + \frac{1}{2}{\rho }^{2} - \Pi,\\
\rho_{t}+(\rho u)_{x}=0.\end{array}\right. \\
\end{equation}
Regularizing this system, we get
\begin{equation}
\left\{\begin{array}{ll}
{u}_{t}^{\varepsilon } + {( u{u}_{x}) }^{\varepsilon } + \Pi_{x}^{\varepsilon } = 0,\\
 {u}_{tx}^{\varepsilon} + {( u   {u}_{xx}) }^{\varepsilon} + \frac{1}{2}{( {u}_{x}^{2}) }^{\varepsilon} = {( {u}^{2}) }^{\varepsilon} + \frac{1}{2}{( {\rho }^{2}) }^{\varepsilon} - {\Pi}^{\varepsilon},\\
{\rho }_{t}^{\varepsilon} + {( \rho u) }_{x}^{\varepsilon} = 0.\end{array}\right. \\
\end{equation}
Multiplying $(3.4)_2$ and $(3.5)_2$ by $u_{x}^\varepsilon$ and $u_{x}$ respectively, we have
\[\left\{  \begin{array}{l} {u}_{xt}   {u}_{x}^{\varepsilon } + {u}{u}_{xx} {u}_{x}^{\varepsilon } + \frac{1}{2}{u}_{x}^{2}   {u}_{x}^{\varepsilon } = ( {{u}^{2} + \frac{1}{2}{\rho }^{2} - \Pi}) {u}_{x}^{\varepsilon }, \\  u^\varepsilon _{tx} {u}_{x} + {( u  {u}_{xx}) }^{\varepsilon}  {u}_{x} + \frac{1}{2}{( {u}_{x}^{2}) }^{\varepsilon}   {u}_{x} = [ {({u}^{2})^\varepsilon + \frac{1}{2}({\rho }^{2})^\varepsilon - \Pi^\varepsilon}] {u}_{x}.\end{array}\right.\]
Multiplying $(3.4)_3$ and $(3.5)_3$ by $\rho^\varepsilon$ and $\rho$ respectively, we have
\[\left\{  \begin{array}{l}   \rho_{t} \rho^\varepsilon+(\rho u)_{x}\rho^ \varepsilon=0 , \\ \rho_{t}^\varepsilon \rho+ (\rho u)_{x}^\varepsilon \rho=0 . \end{array}\right.\]
Adding the above four equations, we obtain
\begin{align}
& ( u_x u_x^{\varepsilon} )_t + ( \rho \rho^{\varepsilon} )_t
- \underbrace{ [u_x^{\varepsilon} u^2+u_x ( u^2 )^{\varepsilon}]}_{\mathrm{I}} \nonumber \\
& + \underbrace{u_x^{\varepsilon} u u_{xx} + \frac{1}{2} u_x^{\varepsilon} u_x^2 + u_x ( u u_{xx} )^{\varepsilon} + \frac{1}{2} u_x( u_x^2 )^{\varepsilon}}_{\mathrm{II}}\nonumber \\
& + \underbrace{\rho^{\varepsilon} ( \rho u )_x + \rho ( \rho u )_x^{\varepsilon}-\frac{1}{2} u_x^{\varepsilon} \rho^2 -u_x ( \frac{1}{2} \rho^2 )^{\varepsilon}}_{\mathrm{III}} \nonumber \\
&=  - u_x^{\varepsilon} \Pi - u_x \Pi^{\varepsilon}.
\end{align}
To rewrite the expression for II and III, we apply product rule identities and algebraic rearrangement to restructure the terms as follows
\begin{equation}
\begin{split}
\mathrm{II} &= u_x^{\varepsilon} u u_{xx}
+ \frac{1}{2} u_x^{\varepsilon}  u_x^2
+ u_x  ( u u_{xx} )_x^{\varepsilon}
+ \frac{1}{2} u_x ( u_x^2 )^{\varepsilon} \\
   &= u_x^{\varepsilon} (u u_x)_x-u_x^{\varepsilon} u^2_x
+ \frac{1}{2} u_x^{\varepsilon}  u_x^2
+ u_x  ( u u_x )_x^{\varepsilon}-u_x (u^2_x)^\varepsilon
+ \frac{1}{2} u_x  ( u_x^2 )^{\varepsilon} \\
&=u_x^{\varepsilon} (u u_x)_x-\frac{1}{2} u_x^{\varepsilon}  u_x^2
+ u_x  ( u u_x )_x^{\varepsilon}- \frac{1}{2} u_x ( u_x^2 )^{\varepsilon}\\
   &= [ u_x^{\varepsilon}  ( u u_x ) ]_x
- u_{xx}^{\varepsilon} ( u u_x )
+ [ u_x  ( u u_x )^{\varepsilon} ]_x \\
   &\quad - u_{xx}  ( u u_x )^{\varepsilon}
- \frac{1}{2} u_x^{\varepsilon}  u_x^2
- \frac{1}{2} u_x  ( u_x^2)^{\varepsilon},
\end{split}
\end{equation}
\begin{equation}
\begin{split}
\mathrm{III} &= \rho^{\varepsilon} ( \rho u )_x
+ \rho ( u \rho )_x^{\varepsilon}
- \frac{1}{2} u_x^{\varepsilon} \rho^2
- u_x  ( \frac{1}{2} \rho^2 )^{\varepsilon} \\
&= [ \rho^{\varepsilon} ( \rho u ) ]_x
- \rho_x^{\varepsilon}  ( \rho u )
+ [ \rho ( u \rho )^{\varepsilon} ]_x
- \rho_x ( u \rho )^{\varepsilon}
- \frac{1}{2} u_x^{\varepsilon}  \rho^2
- \frac{1}{2} u_x ( \rho^2 )^{\varepsilon}.
\end{split}
\end{equation}
Substituting $(4.29)-(4.30)$ into $(4.28)$, we obtain
$${( {u}_{x}{u}_{x}^{\varepsilon})}_{t} + {( \rho{\rho}^{\varepsilon}) }_{t}-u_x^{\varepsilon} u^2-u_x ( u^2 )^{\varepsilon}+ {( {u}_{x}^{\varepsilon }   u{u}_{x}) }_{x} - {u}_{xx}^{\varepsilon }u{u}_{x} + {[ {u}_{x}   {( u{u}_{x}) }^{\varepsilon }] }_{x} $$
$$ - {u}_{xx}{( u{u}_{x}) }^{\varepsilon } - \frac{1}{2}{u}_{x}^{\varepsilon }   {u}_{x}^{2} - \frac{1}{2}{u}_{x}   {( {u}_{x}^{2}) }^{\varepsilon }- {\rho }_{x}^{\varepsilon }   {\rho u}+ {[ {\rho }^{\varepsilon }   ( \rho u) ] }_{x} + {[ \rho {( u\rho ) }^{\varepsilon }] }_{x} - {\rho }_{x}{( u\rho ) }^{\varepsilon }$$
$$ - \frac{1}{2}{u}_{x}^{\varepsilon }   {\rho }^{2} - \frac{1}{2}{u}_{x}  {( {\rho }^{2}) }^{\varepsilon } +u_{x}^{\varepsilon}\Pi+u_{x}\Pi^{\varepsilon}=0,$$
it follows that\\
\begin{equation}
\begin{split}
&(u_x u_x^\varepsilon)_t + (\rho \rho^\varepsilon)_t
+ (u_x^\varepsilon u u_x)_x + [u_x(u u_x)^\varepsilon]_x \\
&+ [\rho^\varepsilon (\rho u)]_x + \left[\rho(u\rho)^\varepsilon\right]_x
+u_{x}^{\varepsilon}\Pi+u_{x}\Pi^{\varepsilon} \\
&= u_x^\varepsilon u^2 + u_x (u^2)^\varepsilon
+ u_{xx}^\varepsilon u u_x + u_{xx}(u u_x)^\varepsilon
+ \frac{1}{2} u_x (u_x^2)^\varepsilon + \frac{1}{2} u_x^\varepsilon u_x^2 \\
&+ \rho_x^\varepsilon \rho u + \rho_x (u\rho)^\varepsilon
+ \frac{1}{2} u_x^\varepsilon \rho^2 + \frac{1}{2} u_x (\rho^2)^\varepsilon.
\end{split}
\end{equation}
For the sake of simplicity, we define $\delta u(l)=u(x+l)-u(x)=f-u$, $\delta u_x(l)=u_x(x+l)-u_x(x)=f _x-u_ x $ and $\delta \rho(l)=\rho(x+l)-\rho(x)=g-\rho$, by use of them we can abtain\\
\begin{equation}
\begin{split}
\int_{\mathbb{S}} \varphi_{\varepsilon}'(t) [ \delta u(l) ]^{3} dl
&= \int_{\mathbb{S}} \varphi_{\varepsilon}'(l) [ u(x+l) - u(x) ]^{3} dl \\
&= \int_{\mathbb{S}} \varphi_{\varepsilon}'(l) [ u^{3}(x+l) - 3u^{2}(x+l)u(x) \\
&\quad + 3u(x+l)u^{2}(x) - u^{3}(x) ] dl.
\end{split}
\end{equation}
By performing the change of variables $l=\gamma-x$ and carrying out direct computation, we have\\
\[
\begin{aligned}
\int_{\mathbb{S}} \varphi_{\varepsilon}'(l) u^{3}(x+l) dl
&= \int_{\mathbb{S}} \partial_{\gamma} \varphi_{\varepsilon}(\gamma - x) u^{3}(\gamma) d\gamma \\
&= - \int_{\mathbb{S}} \partial_{x} \varphi_{\varepsilon}(\gamma - x) u^{3}(\gamma) d\gamma \\
&= - \partial_{x} \int_{\mathbb{S}} \varphi_{\varepsilon}(\gamma - x) u^{3}(\gamma) d\gamma \\
&= - ( \varphi_{\varepsilon} * u^{3} )_{x}
= - ( u^{3} )_{x}^{\varepsilon}.
\end{aligned}
\]
Similarly, we get\\
\[
\begin{aligned}
\int_{\mathbb{S}} \varphi_{\varepsilon}'(l) [ -3u^{2}(x+l)u(x) ] dl
&= -3u(x) \int_{\mathbb{S}} \varphi_{\varepsilon}'(l) u^{2}(x+l) dl \\
&= ( -3u(x) ) [ - ( u^{2} )_{x}^{\varepsilon} ]
= 3u(x) ( u^{2} )_{x}^{\varepsilon},
\end{aligned}
\]

\[
\begin{aligned}
\int_{\mathbb{S}} \varphi_{\varepsilon}'(l) \left[ 3u(x+l)u^{2}(x) \right] dl
&= 3u^{2}(x) \int_{\mathbb{S}} \varphi_{\varepsilon}'(l) u(x+l) dl \\
&= -3u^{2} u^{\varepsilon}_x,
\end{aligned}
\]

\[\int_{\mathbb{S}}{\varphi }_{\varepsilon }^{\prime }\left( l\right) u^3(x) {dl}=0.\]
Substituting these four identities into $(4.32)$, we obtain\\
\begin{equation}
\begin{aligned}
\int_{\mathbb{S}} \varphi_{\varepsilon}'(l) [ \delta u(l) ]^{3} dl
&= - ( u^{3} )_{x}^{\varepsilon}
+ 3u ( u^{2} )_{x}^{\varepsilon}
- 3u^{2} u_{x}^{\varepsilon} \\
&= - ( u^{3} )_{x}^{\varepsilon}
+ 3u ( u^{2} )_{x}^{\varepsilon}
- 3u^{2} u_{x}^{\varepsilon} \\
&= - ( u^{3} )_{x}^{\varepsilon}
+ 3 [u(u^2)^{\varepsilon}]_x
- 3u_{x} ( u^{2} )^{\varepsilon}
- 3u^{2} u_{x}^{\varepsilon},
\end{aligned}
\end{equation}
thus\\
\begin{equation}
u^\varepsilon_x u^2+u_x (u^2)^\varepsilon=-\frac{1}{3} \int _{\mathbb{S}}{\varphi }_{\varepsilon}'( l)[ {\delta u}_{}( {l}) ]^3{dl}-\frac{1}{3}(u^3)^\varepsilon_x+[u({u^2})^\varepsilon)]_x.
\end{equation}
By suitable modification of the deduction of the above processes, we can get\\
\begin{equation}
\begin{aligned}
&\int_{\mathbb{S}} \varphi_{\varepsilon}'(l) \delta u(l)[ \delta u_x(l)]^{2} dl\\
=&\int_{\mathbb{S}} \varphi_{\varepsilon}'(l) [u(x+l)u_x^2(x+l)
- 2u(x+l)u_x(x+l)u_x(x)
+ u(x+l)u_x^2(x) \\
&\quad
- u(x)u_x^2(x+l)
+ 2u(x)u_x(x+l)u_x(x)
- u(x)u_x^2(x)
] dl \\
=& - ( u u_x^2)_{x}^{\varepsilon}+ 2u_x ( u u_x )_{x}^{\varepsilon}
- u_x^2 u_{x}^{\varepsilon}\\
&+ u ( u_x^2 )_{x}^{\varepsilon}- 2u u_x ( u_x )_{x}^{\varepsilon}
+ \int_{\mathbb{S}} \varphi_{\varepsilon}'(l) u u_x^2 dl \\
=& - ( u u_x^2 )_{x}^{\varepsilon}+ 2u_x ( u u_x )_{x}^{\varepsilon}- u_x^2 u_{x}^{\varepsilon}
+ u ( u_x^2)_{x}^{\varepsilon}- 2u u_x ( u_x )_{x}^{\varepsilon}.
\end{aligned}
\end{equation}
Since
\[2{u}_{x}{\left( u{u}_{x}\right) }_{x}^{\varepsilon }+{u} {\left( {u}_{x}^{2}\right) }_{x}^{\varepsilon }=2[(u u_x)^\varepsilon u_x]_x-2(u u_x)^\varepsilon  u_{xx}+[u({u_x^2})^{\varepsilon }]_x-u_x(u_x^2)^\varepsilon,\]
we obtain\\
\begin{equation}
\begin{aligned}
& \int_{\mathbb{S}} \varphi_{\varepsilon}'(l) [ u(x+l) - u(x) ] [ u_x(x+l) - u_x(x) ]^{2} dl \\
= &- ( u u_x^2 )_{x}^{\varepsilon}
+ 2 [ ( u u_x )^{\varepsilon} u_x ]_{x}
- 2 ( u u_x )^{\varepsilon} u_{xx}
+ [ u ( u_x^2 )^{\varepsilon} ]_{x} \\
&- u_x ( u_x^2 )^{\varepsilon}
- u_x^2 u_{x}^{\varepsilon}
- 2 u u_x ( u_x )_{x}^{\varepsilon}.
\end{aligned}
\end{equation}
Rearranging terms, we have\\
\[
\begin{aligned}
& 2( u u_x )^{\varepsilon} u_{xx}
+ u_x ( u_x^2 )^{\varepsilon}
+ u_x^2 u_{x}^{\varepsilon}
+ 2u u_x ( u_x )_{x}^{\varepsilon} \\
=& - \int_{\mathbb{S}} \varphi_{\varepsilon}'(l) [u(x+l) - u(x) ] [ u_x(x+l) - u_x(x) ]^{2} dl \\
&- ( u u_x^2 )_{x}^{\varepsilon}
+ 2[ ( u u_x )^{\varepsilon} u_x ]_{x}
+ [ u( u_x^2 )^{\varepsilon} ]_{x},
\end{aligned}
\]
it follows that\\
\begin{equation}
\begin{aligned}
& ( u u_x )^{\varepsilon} u_{xx}
+ u u_x  (u_{x})_x ^{\varepsilon}
+ \frac{1}{2}  u_x ( u_x^2 )^{\varepsilon} + \frac{1}{2}u_x^2 u_{x}^{\varepsilon}  \\
=& -\frac{1}{2} \int_{\mathbb{S}} \varphi_{\varepsilon}'(l) \delta u(l) [ \delta u_x(l) ]^{2} dl
- \frac{1}{2} ( u u_x^2 )_{x}^{\varepsilon}
+ [ ( u u_x )^{\varepsilon} u_x ]_{x}
+ \frac{1}{2} [ u ( u_x^2 )^{\varepsilon} ]_{x}.
\end{aligned}
\end{equation}
Repeating the above argument, we obtain\\
\begin{equation}
\begin{aligned}
& \int_{\mathbb{S}} \varphi_{\varepsilon}'(l) \delta u(l) [ \delta \rho(l) ]^{2} dl \\
= &\int_{\mathbb{S}} \varphi_{\varepsilon}'(l) ( f - u ) ( g - \rho )^{2} dl \\
= &\int_{\mathbb{S}} \varphi_{\varepsilon}'(l) ( f - u ) ( g^{2} - 2g\rho + \rho^{2} ) dl \\
= &\int_{\mathbb{S}} \varphi_{\varepsilon}'(l) ( fg^{2} - 2fg\rho + f\rho^{2} - ug^{2} + 2g\rho u - u\rho^{2} ) dl \\
= &- ( u\rho^{2} )_{x}^{\varepsilon}
+ 2\rho ( u\rho )_{x}^{\varepsilon}
- \rho^{2} u_{x}^{\varepsilon}
+ u ( \rho^{2} )_{x}^{\varepsilon}
- 2u\rho \rho_{x}^{\varepsilon}.
\end{aligned}
\end{equation}
A routine computation yields\\
\[{2\rho }{( u\rho ) }_{x}^{\varepsilon} + u{( {\rho }^{2}) }_{x}^{\varepsilon }
= 2[\rho(\rho u)^{\varepsilon}]_x - 2{\rho }_{x}  {( u{\rho }) }^{\varepsilon} - {u}_{x}{( {\rho}^{2}) }^\varepsilon+[u(\rho^2)^{\varepsilon}]_{x}.\]
Substituting this expression into $(4.38)$, we obtain\\
\[
\begin{aligned}
& \int_{\mathbb{S}} \varphi_{\varepsilon}'(l) \delta u(l) [ \delta \rho(l) ]^{2} dl \\
=& - ( u\rho^{2} )_{x}^{\varepsilon} + 2 [ \rho ( u\rho )^{\varepsilon} ]_{x} - 2\rho_{x} ( u\rho )^{\varepsilon} + [ u ( \rho^{2} )^{\varepsilon} ]_{x} - u_{x} ( \rho^{2} )^{\varepsilon} - \rho^{2} u_{x}^{\varepsilon} - 2u\rho \rho_{x}^{\varepsilon}.
\end{aligned}
\]
So\\
\begin{equation}
\begin{aligned}
& \rho_{x}^{\varepsilon} \rho u + \rho_{x} ( u\rho )^{\varepsilon} + \frac{1}{2} u_{x}^{\varepsilon} \rho^{2} + \frac{1}{2} u_{x} ( \rho^{2} )^{\varepsilon} \\
=& -\frac{1}{2} \int_{\mathbb{S}} \varphi_{\varepsilon}'(l) \delta u(l) [ \delta \rho(l) ]^{2} dl
- \frac{1}{2} ( u\rho^{2} )_{x}^{\varepsilon}
+ [ \rho ( u\rho )^{\varepsilon} ]_{x}
+ \frac{1}{2} [ u ( \rho^{2} )^{\varepsilon}]_{x}.
\end{aligned}
\end{equation}
Substituting $(4.34)$, $(4.37)$ and $(4.39)$ into $(4.31)$, we have\\
\[
\begin{aligned}
& (u_x u^\varepsilon_x)_t+(\rho \rho^\varepsilon)_t+(u^\varepsilon_x u u_x)_x+[u_x (u u_x)^\varepsilon]_x \\
& +[\rho^\varepsilon(\rho u)]_x+[\rho(u \rho)^\varepsilon]_x +u_{x}^{\varepsilon}\Pi+u_{x}\Pi^{\varepsilon} \\
= &-\frac{1}{3}\int_{\mathbb{S}}\varphi_\varepsilon'(l)\bigl[\delta u(l)\bigr]^3 dl
-\frac{1}{3}(u^3)^\varepsilon_x+[u(u^2)^\varepsilon]_x \\
&-\frac{1}{2}\int_{\mathbb{S}}\varphi_\varepsilon'(l)\delta u(l)\bigl[\delta u_x(l)\bigr]^2 dl
-\frac{1}{2}(u u_x^2)^\varepsilon_x+\bigl[(u u_x)^\varepsilon u_x\bigr]_x
+\frac{1}{2}[u(u_x^2)^\varepsilon]_x \\
&-\frac{1}{2}\int_{\mathbb{S}}\varphi_\varepsilon'(l)\delta u(l)\bigl[\delta \rho(l)\bigr]^2 dl
-\frac{1}{2}(u\rho^2)^\varepsilon_x+\bigl[\rho(u\rho)^\varepsilon\bigr]_x
+\frac{1}{2}[u(\rho^2)^\varepsilon]_x,
\end{aligned}
\]
it follows that\\
\begin{equation}
\begin{aligned}
& (u_{x} u_{x}^{\varepsilon})_{t} + (\rho \rho^{\varepsilon})_{t}
- [ u ( u^{2} )^{\varepsilon} ]_{x}
+ ( u_{x}^{\varepsilon} u u_{x} )_{x} \\
&+[ \rho^{\varepsilon} (\rho u) ]_{x}
+u_{x}^{\varepsilon}\Pi+u_{x}\Pi^{\varepsilon}
+ \frac{1}{3} ( u^{3} )_{x}^{\varepsilon}
+ \frac{1}{2} ( u u_{x}^{2} )_{x}^{\varepsilon}
- \frac{1}{2} [ u ( u_{x}^{2} )^{\varepsilon} ]_{x} \\
&+ \frac{1}{2} ( u \rho^{2} )_{x}^{\varepsilon}
- \frac{1}{2} [ u ( \rho^{2} )^{\varepsilon} ]_{x} \\
=& -\frac{1}{2} \int_{\mathbb{S}} \varphi_{\varepsilon}'(l) [ \delta u(l) ]^{3} dl
- \frac{1}{2} \int_{\mathbb{S}} \varphi_{\varepsilon}'(l) \delta u(l) [ \delta u_{x}(l) ]^{2} dl \\
&\quad - \frac{1}{2} \int_{\mathbb{S}} \varphi_{\varepsilon}'(l) \delta u(l) [ \delta \rho(l) ]^{2} dl
\end{aligned}
\end{equation}
Letting $\varepsilon\rightarrow 0$ on both sides of $(4.40)$, and noting $(4.15)$, $(4.20)$, $and (4.21)$, we obtain
\[(u_{x}^{2}+\rho^{2})_{t}+[u(u_{x}^{2}+\rho^{2})]_{x}=\frac{2}{3}(u^{3})_{x}-2\Pi u_{x}=2(u^{2}-\Pi)u_{x}.\]
This shows that the weak solution $(u,\rho)$ is a conservative solution of system $(2.1)$. By Theorem $2.2$ in \cite{LZ}, this solution is unique.
This completes the proof of Theorem $2.2$.

\end{proofof}

\bigskip
\noindent\textbf{Acknowledgments} This work was partially supported by the National Natural Science Foundation of China
(Nos. 11701525 and 11971446). The
authors thank two referees for their valuable comments and
suggestions.

\end{document}